\documentclass{article}

\usepackage{amsthm,amsmath,amssymb}

\newtheorem{theorem}{Theorem}

\begin{document}

\title{A short and elementary proof of Brauer's theorem}
\markright{Brauer's theorem}
\author{Judith J.~McDonald and Pietro Paparella}

\maketitle

\begin{abstract}
A short and elementary proof is given of a celebrated eigenvalue-perturbation result due to Alfred Brauer.

\emph{MathEduc Subject Classification}: H65


\end{abstract}

In 1952, Alfred Brauer \cite[Theorem 27]{b1952} proved the following useful eigenvalue perturbation result, which is used in, e.g., \emph{deflation techniques} when computing eigenvalues (see, e.g., Saad \cite[Section 4.2]{s2011}) and the longstanding \emph{nonnegative inverse eigenvalue problem} (see, e.g., Julio and Soto \cite{js2020} and references therein).  

\begin{theorem}
Let $A$ be an $n$-by-$n$ matrix with complex entries and suppose that $A$ has eigenvalues \( \{ \lambda, \lambda_2, \dots, \lambda_n \} \) (including multiplicities). If $x$ is an eigenvector associated with $\lambda$ and $y \in \mathbb{C}^n$, then the matrix $A + xy^*$ has eigenvalues $\{ \lambda + y^*x, \lambda_2, \dots, \lambda_n \}$. 
\end{theorem}

Brauer's proof and a recent proof by Melman \cite{m2020} rely on the \emph{principle of biorthogonality}. While both proofs are elementary, they are somewhat lengthy. Meanwhile, Horn and Johnson give two brief proofs, but the first relies on the \emph{adjoint of a matrix} \cite[p.~51]{hj2013} and the second requires a unitary matrix \cite[p.~122]{hj2013} as opposed to just an invertible matrix (see details below). Another proof by Reams \cite[p.~368]{r1996} implicitly relies on \emph{Schur triangularization}. 

The proof below relies on the elementary fact that similar matrices are cospectral, an elementary fact that is covered in a first-course in linear algebra.

\begin{proof}
Let $Q$ be any invertible matrix whose first column is $x$---say $Q = \begin{bmatrix} x & R \end{bmatrix}$. Notice that $Q^{-1} x = e_1$---the first canonical basis vector of $\mathbb{C}$---and that $Q^{-1} A Q e_1 = \lambda e_1$. A simple calculation reveals that 
\begin{align*}
    Q^{-1} (A + xy^*) Q 
    &= Q^{-1} A Q + Q^{-1} xy^* Q \\
    &= \begin{bmatrix} \lambda & u^* \\ 0 & C \end{bmatrix} + e_1y^* Q = \begin{bmatrix} \lambda + y^*x & u^* + y^* R \\ 0 & C \end{bmatrix}
\end{align*}
i.e., $\sigma(A) = \{ \lambda \} \cup \sigma(C)$ and $\sigma(A + xy^*) = \{ \lambda + y^*x \} \cup \sigma(C)$.
\end{proof}

\end{document}